\documentclass{amsproc}

\usepackage{latexsym}
\usepackage{amsmath}
\usepackage{amsfonts}
\usepackage{amssymb}
\usepackage{amscd}

\numberwithin{equation}{section}

\newtheorem{proposition}{Proposition}[section]
\newtheorem{lemma}[proposition]{Lemma}
\newtheorem{theorem}[proposition]{Theorem}
\newtheorem{corollary}[proposition]{Corollary}
\newtheorem{definition}[proposition]{Definition}

\newtheorem{question}[proposition]{Question}

\def\itheorem#1#2{\newtheorem{#1}[proposition]{#2}}

\renewcommand{\labelenumi}{(\roman{enumi})}

\itheorem{Remark}{Remark}

\def\Mod{\mathop{\rm Mod}\nolimits}
\def\Hom{\mathop{\rm Hom}\nolimits}
\def\Ext{\mathop{\rm Ext}\nolimits}

\def\Trans{\mathop{\rm Trans}\nolimits}
\def\seq{\mathop{\rm seq}\nolimits}

\def\rk{\mathop{\rm rk}\nolimits}
\def\dom{\mathop{\rm dom}\nolimits}
\def\rg{\mathop{\rm rg}\nolimits}
\def\l{\mathop{\rm l}\nolimits}
\def\cf{\mathop{\rm cf}\nolimits}

\def\Ord{\mathop{\rm ORD}\nolimits}

\def\dQ{{\overset{.}{Q}}}
\def\dA{{\overset{.}{A}}}
\def\df{{\overset{.}{f}}}
\def\dS{{\overset{.}{S}}}

\def\Q{{\mathbb{Q}}}
\def\Z{{\mathbb{Z}}}

\def\C{{\mathfrak{C}}}

\def\tf{{\mathcal{T}\mathcal{F}}}

\def\P{{\mathcal{P}}}

\begin{document}

\title{Kulikov's problem on universal torsion-free abelian groups}

\author{Saharon Shelah}

\address{The Hebrew University, Givat Ram, Jerusalem 91904,
Israel and \\ Rutgers University, Newbrunswick, NJ U.S.A.}

\email{Shelah@math.huji.ac.il}

\thanks{Publication 772 in the first author's list of publication. The first author was supported
by project No. G-0545-173,06/97 of the {\em German-Israeli
Foundation for Scientific Research \& Development.}}

\author{Lutz Str\"ungmann}

\address{Fachbereich 6 -- Mathematik,
University of Essen, 45117 Essen, Germany}

\curraddr{The Hebrew University, Givat Ram, Jerusalem 91904,
Israel}

\email{lutz@math.huji.ac.il}

\thanks{The second author was supported by a MINERVA fellowship.\\
2000 Mathematics Subject Classification 20K10, 20K20, 20K35}


\begin{abstract}
Let $T$ be an abelian group and $\lambda$ an uncountable regular
cardinal. We consider the question of whether there is a
$\lambda$-universal group $G^*$ among all torsion-free abelian
groups $G$ of cardinality less than or equal to $\lambda$
satisfying $\Ext(G,T)=0$. Here $G^*$ is said to be
$\lambda$-universal for $T$ if, whenever a torsion-free abelian
group $G$ of cardinality less than or equal to $\lambda$
satisfies $\Ext(G,T)=0$, then there is an embedding of $G$ into
$G^*$. For large classes of abelian groups $T$ and cardinals
$\lambda$ it is shown that the answer is consistently no. In
particular, for $T$ torsion, this solves a problem of Kulikov.
\end{abstract}

\maketitle

\section{Introduction}
Given a class $\C$ of objects and a property $P$ it is natural to ask for universal objects in $\C$ with respect to $P$. A universal object with respect to $P$ is an object $C \in \C$ satisfying $P$ such that every other object of the class $\C$ that satisfies $P$ can be identified with a subobject of $C$. The existence of universal objects clearly simplifies the theorie of the objects satisfying $P$ since very often properties of objects are inherited by subobjects. Thus, there is a distinguished object in $\C$ with respect to $P$ that is a representative among all the objects of $\C$ satisfying $P$.\\ On the other hand, if there exist no universal objects with respect to $P$, this indicates that the objects satisfying $P$ have a complicated structure. Since the definition of universal objects is ``universal'' the question of their existence appears in almost every field of mathematics. e.g. group theory, theory of ordered structures, Banach-spaces etc.\\
In the present paper we focus on the theory of abelian groups and define our class $\C=\tf_{\lambda}$ to be the class of all torsion-free abelian groups $G$ of cardinality (rank) less than or equal to $\lambda$, where $\lambda$ is a fixed cardinal. The property $P=P_T$ is related to a fixed torsion abelian group $T$ and a group $G \in \tf_{\lambda}$ satisfies $P_T$ if and only if $\Ext(G,T)=0$, where $\Ext(-,T)$ denotes the first derived functor of the functor $\Hom(-,T)$. It was Kulikov \cite[Question 1.66]{KN} who first asked whether or not there exist universal groups in $\tf_{\lambda}$ for all (uncountable) cardinals $\lambda$ and torsion abelian groups $T$. Clearly, if the group $T$ is cotorsion, hence satisfies $\Ext(\Q,T)=0$, then there is always (for every $\lambda$) a universal group in $\tf_{\lambda}$, namely the torsion-free divisible group of rank $\lambda$. Moreover, since every free abelian group $F$ satisfies $\Ext(F,T)=0$, one has to consider torsion-free groups of cardinality less than or equal to a fixed cardinal rather than searching for universal objects among all torsion-free abelian groups with $P_T$. As was mentioned earlier, the existence of universal objects in $\tf_{\lambda}$ with respect to $P_T$ sheds some light on the complexity of the structure of objects in $\tf_{\lambda}$. This has also consequences to more complicated theories. We shall discuss cotorsion theories as an example.\\
 Cotorsion theories for abelian groups have been introduced by
Salce in 1979 \cite{S}. Following his notation we call a pair
$(\mathcal{F}, \mathcal{C})$ a cotorsion theory if $\mathcal{F}$
and $\mathcal{C}$ are classes of abelian groups which are maximal
with respect to the property
that $\Ext(F,C)=0$ for all $F \in \mathcal{F}$, $C \in \mathcal{C}$.\\
Salce \cite{S} has shown that every cotorsion theory is
cogenerated by a class of torsion and torsion--free groups where
$(\mathcal{F},\mathcal{C})$ is said to be cogenerated by the
class $\mathcal{A}$ if $\mathcal{C}=\mathcal{A}^{\perp}=\{ X\in
\Mod$--$\Z \mid \Ext(A,X)=0 \textit{ for all } A \in \mathcal{A}
\}$ and $\mathcal{F} = ^{\perp}( \mathcal{A}^{\perp})=\{ Y \in
\Mod$--$\Z \mid \Ext(Y,X)=0 \textit{ for all } X \in
\mathcal{A}^{\perp} \}$. Examples for cotorsion theories are:
\\ $(\mathcal{L},
\Mod$--$\Z)=(^{\perp}(\Z^{\perp}), \Z^{\perp})$ where
$\mathcal{L}$ is the class of all free groups, and the classical
one $(\mathcal{TF}, \mathcal{CO})=(^{\perp}(\Q^{\perp}),
\Q^{\perp})$ where $\mathcal{TF}$ is the class of all
torsion--free groups and $\mathcal{CO}$ is the class of all
(classical) cotorsion groups. In view of the last example the
classes $\mathcal{F}$ and $\mathcal{C}$ of a cotorsion theory
$(\mathcal{F},\mathcal{C})$ are said to be the torsion--free
class and
the cotorsion class of this cotorsion theory. \\
If we restrict our attention to cotorsion classes
cogenerated by a single torsion-free group $G$, then ordering these classes by
inclusion, we obviously have that $\Z^{\perp}$ is maximal and $\Q^{\perp}$ is minimal among these classes. Moreover, in \cite{GSW} G\"obel, Wallutis and the first author have shown that any partially ordered set can be embedded into the lattice of all cotorsion classes. Hence there is no hope at all to characterize these classes. But if we restrict to torsion-free groups in $\tf_{\lambda}$, then the existence of a universal group can be helpful. If $G \in \tf_{\lambda}$ and $T \in G^{\perp}$ is ``as complicated as possible'', then the universal object $C$ related to $P_T$ satisfies $C^{\perp} \subseteq G^{\perp}$ and hence we obtain new information about $G^{\perp}$. This is just one example where universal objects could be helpful.\\
To the authors's knowledge there is no published literature on Kulikov's problem except for \cite{St}, where the second author proved that in G\"odel's universe ($V=L$) for every cardinal $\lambda$ and torsion abelian group $T$ there exists a universal group $G \in \tf_{\lambda}$ if $T$ has only finitely many non-trivial bounded $p$-components. Moreover, if $\lambda$ is finite, then this is an if and only if result already in ZFC. Thus, in $V=L$, the class $\tf_{\lambda}$ behaves well with respect to the property $T_P$ for a large class of torsion abelian groups $T$.\\
In this paper we prove that the result from \cite{St} is not provable in ZFC. We show that it is consistent with ZFC and $GCH$ that for every abelian group $T$ (not necessarily torsion) and every uncountable regular cardinal $\lambda$ there is a cardinal $\kappa > \lambda$ such that the class $\tf_{\kappa}$ has no universal object with respect to $P_T$. Moreover, we prove that for torsion abelian groups $T$ of cardinality less than or equal to $\aleph_1$ (also more general situations are considered) there is no uncountable cardinal $\lambda$ such that $\tf_{\lambda}$ has universal groups with respect to the property $P_T$. This answers Kulikov's problem consistently in the negative.\\

All groups under consideration are abelian. The notations are
standard and for unexplained notions in abelian group theory and
set theory we refer to \cite{Fu} and \cite{EM}, \cite{J}. For
uniformization see \cite{K} or \cite{S4}.

\section{$\lambda$-universal groups}
In this section we introduce the notions of $\lambda$-universal
groups for a given group $T$ and obtain some basic properties. We
are mainly interested in the case when our group $T$ is torsion
but for the sake of generality we let $T$ be an arbitrary (abelian) group.
By $\tf$ we denote the class of all torsion-free groups. For a
cardinal $\lambda$ we denote by $\tf_{\lambda}(T)$ the
class of all torsion-free groups $G$ of rank less than or equal
to $\lambda$ such that $\Ext(G,T)=0$, i.e. $\tf_{\lambda}(T)=\{ G
\in \tf : \Ext(G,T)=0 \text{ and } \rk(G) \leq \lambda \}$.
Moreover, we let
$\tf(T)=\bigcup\limits_{\lambda}\tf_{\lambda}(T)=\{ G \in \tf :
\Ext(G,T)=0 \}$ be the class of all torsion-free groups $G$
satisfying $\Ext(G,T)=0$.

\begin{definition}
\label{universaldef}
Let $T$ be a group and $\lambda$ a cardinal. A
torsion-free group $G$ of rank less than or equal
to $\lambda$ is called {\rm $\lambda$-universal for $T$} if $G
\in \tf_{\lambda}(T)$ and for every $H \in \tf_{\lambda}(T)$
there is an embedding $i: H \longrightarrow G$ of $H$ into $G$.
\end{definition}

Note that in Definition \ref{universaldef} for $\lambda$ infinite we may replace the rank of
$G$ by its cardinality. Kulikov \cite[Question 1.66]{KN} asked the
following question:

\begin{question}
Let $T$ be a torsion group and $\lambda$ an uncountable
cardinal. Is there always a $\lambda$-universal group for $T$?
\end{question}

Let us first mention that a positive consistency result was
obtained by the second author in \cite{St}. Moreover, the case of
finite $\lambda$ was considered.

\begin{lemma}[\cite{St}]
\label{fruniversal} Let $T$ be a torsion group and $n$ a strictly
positive integer. Then there exists an $n$-universal group $G$ for
$T$ if and only if $T$ has only finitely many non-trivial bounded
$p$-components. In this case, $G$ is completely decomposable.
\end{lemma}

\begin{lemma}[\cite{St}, $V=L$]
If $T$ is a torsion group with only finitely many non-trivial
bounded $p$-components and $\lambda$ is a cardinal, then there is
a $\lambda$-universal group for $T$ which is completely
decomposable.
\end{lemma}

We therefore shall restrict ourselves to uncountable (regular)
cardinals in most of the results. Let us begin with some basic
observations and recall that a basic subgroup of a torsion group $T$ is a direct sum $B$ of cyclic groups which is pure in $T$ and has divisible quotient $T/B$.

\begin{lemma}
Let $T$ be a torsion group, $B$ a basic subgroup of $T$ and $\lambda$ a cardinal. A torsion-free group $G$ is $\lambda$-universal for $T$ if and only if $G$ is  $\lambda$-universal for $B$.
\end{lemma}

\begin{proof}
The proof follows easily since for any torsion-free group $H$ we
have $\Ext(H,T)=0$ if and only if $\Ext(H,B)=0$ (see for example
\cite[Lemma 1.2]{St}).
\end{proof}

\begin{lemma}
Let $T$ be a group and $T=D \oplus R$ where $D$ is divisible and $R$ is reduced. If $\lambda$ is a cardinal, then a torsion-free group $G$ is $\lambda$-universal for $T$ if and only if $G$ is  $\lambda$-universal for $R$.
\end{lemma}

\begin{proof} The proof is straightforward since for any torsion-free
group $H$ we have $\Ext(H,T)=\Ext(H,D) \oplus \Ext(H,R)$ and
$\Ext(H,D)=0$ since $D$ is divisible. \end{proof}

Thus it is enough to consider reduced groups. Moreover, among the
reduced ones we only have to deal with groups that are not
cotorsion. Recall that a group $T$ is called {\it cotorsion} if
$\Ext(\Q,T)=0$ which is equivalent to $\tf(T)=\tf$.

\begin{lemma}
Let $T$ be a cotorsion group. Then there is a $\lambda$-universal
group for every cardinal $\lambda$.
\end{lemma}

\begin{proof} If $T$ is cotorsion, then for every cardinal $\lambda$ we
have $\bigoplus\limits_{\lambda}\Q \in \tf_{\lambda}(T)$. Since
every torsion-free group of rank less than or equal to
$\lambda$ can be embedded into its divisible hull it follows that
$\bigoplus\limits_{\lambda}\Q $ is $\lambda$-universal for $T$.
\end{proof}

The following lemma shows that it makes sense to restrict
ourselves to (torsion) groups $T$ and cardinals $\lambda$ such
that $\lambda \geq |T|$.

\begin{lemma}[\cite{St}]
\label{pure} Let $G$ be any group and $T$ a torsion group. Then
$\Ext(G,T)=0$ if and only if $\Ext(G,T^{\prime})=0$ for all pure
subgroups $T^{\prime}$ of $T$ such that $|T^{\prime} |
\leq |G|$.
\end{lemma}

We shall even assume that $\lambda > | T |$.

\section{$(T,\lambda,\gamma)$-suitable groups}

In what follows let $\lambda > \gamma$ be fixed infinite regular
cardinals unless otherwise stated.

\begin{definition}
\label{defsuitable} Let $T$ be a group of cardinality less than
$\lambda$. A group $G$ is called $(T,\lambda,\gamma)$-suitable if
the following conditions are satisfied:
\begin{enumerate}
\item $\Ext(G,T)\not=0$;
\item There are free groups $F, F_i$ ($i \leq \gamma$) such that
\begin{enumerate}
\item the $F_i$'s ($i \leq \gamma$) form an increasing chain such
that $F_{\gamma}=\bigcup\limits_{i < \gamma}F_i \subseteq F$;
\item $\rk(F_i) \leq | i | + \aleph_0$;
\item $F/F_i$ is free for all $i < \gamma$;
\item $F/F_{\gamma} \cong G$.
\end{enumerate}
\end{enumerate}
\end{definition}

Our first lemma shows that there is always a
$(T,\lambda,\omega)$-suitable group for non-trivial (not
cotorsion) $T$.

\begin{lemma}
\label{notcotorsion} Let $T$ be a group of cardinality less than
$\lambda$ and $G$ a countable group such that $\Ext(G,T)\not=0$.
Then $G$ is $(T,\lambda,\omega)$-suitable. In particular, if $T$
is not cotorsion then there is a $(T,\lambda,\gamma)$-suitable
group $G$.
\end{lemma}

\begin{proof} Let $T$ and $G$ be as stated. Choose a free resolution \[0
\rightarrow K \overset{id}{\rightarrow} F \rightarrow G
\rightarrow 0 \] of $G$. Without loss of generality we may assume
that $K$ and $F$ are of countable rank. Choose elements $e_i \in
K$ ($i < \omega$) such that $K=\bigoplus\limits_{i < \omega}\Z
e_i$. Put $F_n=\bigoplus\limits_{i \leq n}\Z e_i \subseteq F$ for $n < \omega$.
Then each $F_n$ is a direct summand of $F$ and hence $G$ is
$(T,\lambda,\omega)$-suitable. If $T$ is not cotorsion, then
$\Ext(\Q,T)\not=0$ and hence the above arguments show that $\Q$ is
$(T,\lambda,\omega)$-suitable. \end{proof}

The next results show the existence of
$(T,\lambda,\gamma$)-suitable groups for uncountable $\gamma$
under certain assumptions. Recall that a group $G$ is called {\it
almost-free} if all its subgroups of smaller cardinality are free.

\begin{lemma}
Let $T$ be a group of cardinality less than $\lambda$
and $G$ an almost-free group of cardinality $\gamma$ such that
$\Ext(G,T)\not=0$. Then $G$ is $(T,\lambda,\gamma)$-suitable.
\end{lemma}

\begin{proof}
Take a $\gamma$-filtration $G=\bigcup\limits_{\alpha < \gamma} G_{\alpha}$ of $G$ such that each $G_{\alpha}$ is free. By \cite[Lemma XII.1.4]{EM} there is a free resolution associated with this filtration. This is to say there are free groups $K=\bigoplus\limits_{\alpha < \gamma}F_{\alpha}$ and $K=\bigoplus\limits_{\alpha < \gamma}K_{\alpha}$ such that the short sequences
\[ 0 \longrightarrow K \longrightarrow F \longrightarrow G \longrightarrow 0 \]
and
\[ 0 \longrightarrow \bigoplus\limits_{\alpha < \beta}K_{\alpha} \longrightarrow \bigoplus\limits_{\alpha < \beta}F_{\alpha} \longrightarrow G_{\beta} \longrightarrow 0 \]
are exact for all $\beta < \gamma$. Since each $G_{\beta}$ is free
it follows that $G$ is $(T,\lambda,\gamma)$-suitable. \end{proof}

\begin{lemma}
\label{epi}
Let $T$ be a group and $H$ an epimorphic image of $T$.
If $G$ is $(H,\lambda,\gamma)$-suitable then $G$ is
$(T,\lambda,\gamma)$-suitable.
\end{lemma}

\begin{proof} The claim follows immediately noting that $\Ext(G,T)=0$
implies $\Ext(G,H)=0$. \end{proof}

\begin{proposition}
\label{diamondsuitable} Let $S\subseteq \gamma$ be stationary
non-reflecting such that $\cf(\alpha)=\omega$ for all $\alpha \in
S$ and assume that $\Diamond_S$ holds. Let $T$ be a group which
is not cotorsion and has an epimorphic image of size less than or
equal to $\gamma$ that is not cotorsion. Then there exists a
strongly $\gamma$-free torsion-free group $G$ of size $\gamma$
which is $(T,\lambda,\gamma)$-suitable. In particular, this holds
if $T$ is torsion or $T$ itself is of cardinality less than or
equal to $\gamma$.
\end{proposition}

\begin{proof} Let $B$ be the epimorphic image of $T$ of size less than
or equal to $\gamma$ which is not cotorsion. Then Lemma \ref{epi}
shows that it is enough to construct a
$(B,\lambda,\gamma)$-suitable group. Therefore, we may assume
without loss of generality that $T$ has cardinality less than or equal to $\gamma$.
Since $T$ is not cotorsion there exists by Lemma
\ref{notcotorsion} a countable torsion-free group $R$ which is
$(T,\lambda,\omega)$-suitable. Let $\lambda_n \geq \aleph_0$ $(n
\in \omega)$ be cardinals. As in \cite[Corollary VII.1.2]{EM}
there exist free abelian groups $K \subseteq F$ such that
$K=\bigcup\limits_{n \in \omega}K_n$ and $F/K_n$ is free for all
$n \in \omega$, $K_0$ is free of rank $\lambda_0$ and
$K_{n+1}/K_n$ is free of rank $\lambda_{n+1}$. Moreover,
$\Ext(F/K,T)\not=0$ since $F/K$ is isomorphic to $R$. As in the
proof of \cite[VII.1.4]{EM} we can construct a torsion-free group
$G$ of cardinality $\gamma$ which has a $\gamma$-filtration
$G=\bigcup\limits_{\alpha < \gamma}G_{\alpha}$ satisfying the
following for all $\alpha < \beta < \gamma$:
\begin{enumerate}
\item $G_{\alpha}$ is free of rank $| \alpha | + \aleph_0$;
\item if $\alpha$ is a limit ordinal, then $G_{\alpha}=\bigcup\limits_{\delta < \alpha}G_{\delta}$;
\item if $\alpha \not\in S$, then $G_{\beta}/G_{\alpha}$ is free of rank $| \beta | + \aleph_0$;
\item if $\alpha \in S$, then $\Ext(G_{\beta}/G_{\alpha},T)\not=0$.
\end{enumerate}
Since $\Diamond_S$ holds it follows that $\Ext(G,T)\not=0$ (see
e.g. \cite[XII.1.15]{EM}), hence $G$ is
$(T,\lambda,\gamma)$-suitable. Finally, if $T$ is torsion then we
choose a basic subgroup $B^{\prime}$ of $T$ and a countable
unbounded direct summand $B$ of $B^{\prime}$ which is therefore
not cotorsion. Note that $B$ exists since $T$ is not cotorsion.
It is well-known that $B$ is an epimorphic image of $T$ (see
\cite[Theorem 36.1]{Fu}), hence Lemma \ref{epi} shows that it is
enough to construct a $(B,\lambda,\gamma)$-suitable group.
\end{proof}

\section{The uniformization}

\relax From now on let $S$ be a stationary subset of $\lambda$ consisting
of limit ordinals of cofinality $\gamma$. To prove our next
theorem we shall use a construction for modules which was almost
identically developed in \cite{ES}. Thus we shall not give all the
proofs but for the convenience of the reader we shall recall the
basic definitions as well as the construction and the main
properties of the constructed module (group).

\begin{definition}
A {\rm ladder system} $\bar{\eta}$ on $S$ is a family of functions
$\bar{\eta}=\left< \eta_{\delta} : \delta \in S \right>$ such that
$\eta_{\delta}: \gamma \rightarrow \delta$ is strictly increasing
with $\sup(\rg(\eta_{\delta}))=\delta$, where $\rg(\eta_{\delta})$
denotes the range of $\eta_{\delta}$. We call the ladder system
{\rm tree-like} if for all $\delta, \nu \in S$ and every $\alpha,
\beta \in \gamma$, $\eta_{\delta}(\alpha)=\eta_{\nu}(\beta)$
implies $\alpha=\beta$ and $\eta_{\delta}(\rho)=\eta_{\nu}(\rho)$
for all $\rho \leq \alpha$.
\end{definition}

For a ladder system $\bar{\eta}=\left< \eta_{\delta} : \delta \in S
\right>$ on $S$ we can form a {\it tree} $B_{\bar{\eta}} \subseteq
{^{\leq{\gamma}}\lambda}$ of height $\gamma$ in the following way:
Let $B_{\bar{\eta}}= \{ \eta_{\delta}\restriction_{\alpha} : \delta \in S,
\alpha \leq \l(\eta_{\delta}) \}$, where $\l(\eta_{\delta})$ denotes
the length of $\eta_{\delta}$, i.e. $l(\eta_{\delta})=\sup\{
\eta_{\delta}(\alpha) : \alpha \in \dom(\eta_{\delta})\}$. Note
that $B_{\bar{\eta}}$ is partially ordered by defining $\eta \leq \nu$ if and
only if $\eta=\nu\restriction_{\l(\eta)}$.
\\
\relax From this tree we now build a group. Let $T$ be a group and let
$G$ be $(T,\lambda,\gamma)$-suitable. Fix a chain
$\left<F_{\alpha} : \alpha \leq \gamma \right>$ for $G$ as in Definition
\ref{defsuitable}. For each $\eta \in B_{\bar{\eta}}$ we let
$H_{\eta}=F_{\l(\eta)}$ and if $\eta \leq \nu \in B_{\bar{\eta}}$ then let
$i_{\eta,\nu}$ be the inclusion map of $H_{\eta}$ into $H_{\nu}$.
Finally, let $H_{\bar{\eta}}$ be the direct limit of $(H_{\eta}, i_{\eta,\nu}
: \eta \leq \nu \in B_{\bar{\eta}})$. More precisely, $H_{\bar{\eta}}$ equals
$\bigoplus\{H_{\eta} : \eta \in B_{\bar{\eta}} \}/K$ where $K$ is the
subgroup generated by all elements of the form $x_{\eta}-y_{\nu}$
where $y_{\nu} \in H_{\nu}, x_{\eta} \in H_{\eta}$, $\eta \leq
\nu$ and $i_{\eta,\nu}(y_{\nu})=x_{\eta}$. Canonically we can
embed $H_{\eta}$ into $H_{\bar{\eta}}$ and we shall therefore regard
$H_{\eta}$ as a submodule of $H_{\bar{\eta}}$ in the sequel.

\begin{definition}
Let $\kappa$ be an uncountable regular cardinal. The tree $B_{\bar{\eta}}$ is
called {\rm $\kappa$-free} if for every $F \subseteq S$ such that
$| F | < \kappa$ there is a function $\Psi: F \rightarrow
\gamma$ such that \[ \{ \{ \eta_{\delta}\restriction_{\alpha} :
\Psi(\delta) < \alpha \leq \gamma \} : \delta \in F \} \] is a
family of pairwise disjoint sets. The ladder system $\bar{\eta}$ is called
{\rm $\kappa$-free} if $B_{\bar{\eta}}$ is $\kappa$-free.
\end{definition}

We now state some properties of the constructed group $H_{\bar{\eta}}$.

\begin{lemma}
\label{lem1} Let $\kappa$ be an uncountable regular cardinal. If
$B_{\bar{\eta}}$ is $\kappa$-free, then $H_{\bar{\eta}}$ is a $\kappa$-free group.
\end{lemma}

\begin{proof} See \cite[Lemma 1.4]{ES}. \end{proof}

\begin{lemma}
\label{lem2} Suppose that $S$ is non-reflecting. Then $H_{\bar{\eta}}$ is not
free.
\end{lemma}

\begin{proof} See \cite[Lemma 1.5]{ES}. Since $S$ is non-reflecting $\bar{\eta}$
is $\lambda$-free and it is easy to see that for all $\delta \in
S$ there exists $\nu \geq \delta$ such that for all $\mu <
\gamma$, $\eta_{\nu}\restriction_{\mu} \in \{
\eta_{\alpha}\restriction_{\mu} : \alpha < \delta \}$.
\end{proof}

We recall the definition of $\mu$-uniformization for a ladder
system $\bar{\eta}$ and a cardinal $\mu$.

\begin{definition}
If $\mu$ is a cardinal and $\bar{\eta}$ is a ladder system on $S$ we say
that $\bar{\eta}$ has {\rm $\mu$-uniformization} if for every family $\{
c_{\delta} : \delta \in S\}$, where $c_{\delta}:
\rg(\eta_{\delta}) \rightarrow \mu$, there exists $\Psi: \lambda
\rightarrow \mu$ and $\Psi^*: S \rightarrow \mu$ such that for
all $\delta \in S$,
$\Psi(\eta_{\delta}(\alpha))=c_{\delta}(\eta_{\delta}(\alpha))$
whenever $\Psi^*(\delta)\leq \alpha < \gamma$.
\end{definition}

\begin{theorem}
\label{uniformization} Let $T$ be a group of cardinality less
than $\lambda$ and $G$ a group which is
$(T,\lambda,\gamma)$-suitable. Moreover, assume that $S$ is
non-reflecting and $\bar{\eta}$ is a tree-like ladder system on $S$ that
has $2^{(| T |^{\gamma})}$-uniformization. Then there
exists a torsion-free group $H$ of size $\lambda$ such that
\begin{enumerate}
\item $H$ has a $\lambda$-filtration $\left< \bar{H}_\alpha : \alpha < \lambda
\right>$;
\item if $\alpha \in S$, then $\bar{H}_{\alpha + 1}/\bar{H}_{\alpha} \cong G$;
\item if $\alpha \not\in S$, then $\bar{H}_{\beta}/\bar{H}_{\alpha}$ is free
for all $\alpha \leq \beta$;
\item $\Ext(H,T)=0$.
\end{enumerate}
In fact, $H$ satisfies $\Ext(H,W)=0$ for all groups $W$ of
cardinality less than or equal to $| T |$.
\end{theorem}

\begin{proof} Let $T$, $S$, $\bar{\eta}$ and $G$ be as stated and choose the
group $H_{\bar{\eta}}=\bigoplus\{ H_{\eta} : \eta \in
B_{\bar{\eta}}\}/K$ as constructed above. Note that $\bar{\eta}$
is $\lambda$-free since $S$ is non-reflecting. Then
$H_{\bar{\eta}}$ is almost-free but not free by Lemma \ref{lem1}
and Lemma \ref{lem2}. Moreover, \cite[Theorem 1.7]{ES} shows that
$H_{\bar{\eta}}$ satisfies $\Ext(H_{\bar{\eta}},W)=0$ for every
group $W$ of size less than or equal to $| T |$. Finally,
it is easy to see that $H_{\bar{\eta}}$ has a
$\lambda$-filtration as stated letting
$\bar{H}_{\alpha}=\bigcup\{(H_{\eta}+K)/K : \eta \in
B_{\bar{\eta}}, \sup(\rg(\eta))<\alpha \}$. \end{proof}

It is our aim to apply Theorem \ref{uniformization} to ladder
systems which have $\mu$-uniformization for all $\mu < \lambda$.
In this case Theorem \ref{uniformization} is applicable to a lot
of regular cardinals $\lambda
> \gamma$ if $T$ is small in some sense (e.g. $| T | \leq
\gamma$). For instance, if $\lambda$ is strongly inaccessible or
$\lambda=\kappa^+$ with $\kappa$ singular and
$\cf(\kappa)=\gamma$, then certainly $\lambda > 2^{| T
|^{\gamma}}$ (even $\kappa > 2^{(|T|^{\gamma})}$) and hence uniformization holds. But problems
arise when $\lambda=\kappa^+$ and $\kappa$ is regular. We shall
show that we can improve Theorem \ref{uniformization} in this
case if $\Diamond_{\gamma}$ holds and the ladder system is
tree-like.

\begin{definition}
\label{defuni} Let $\lambda=\kappa^+$ and $\cf(\kappa)=\kappa$. A
ladder system $\bar{\eta}=\left< \eta_{\delta} : \delta \in S \right>$ has
{\rm strong $\kappa$-uniformization} if for every system $\bar{P}=\left<
P_{\alpha} : \alpha < \lambda \right>$ such that
\begin{enumerate}
\item $\emptyset \not= P_{\delta} \subseteq \{ f | f:
\rg(\eta_{\delta}) \rightarrow \delta \cap \kappa \}$ if $\delta
\in S$;
\item if $\alpha=\eta_{\delta}(i)$ for some $\delta \in S$ and $i
< \gamma$, then $P_{\eta_{\delta}(i)}=\{
f\restriction_{\rg(\eta_{\delta}\restriction_{(i+1)})} | f \in
P_{\delta} \}$;
\item if $\delta \in S$ and $i < \gamma$ is a limit ordinal, then
for every increasing sequence $\left< f_j : j < i \right>$, $f_j
\in P_{\eta_{\delta}(j)}$ there exists $f_i \in
P_{\eta_{\delta}(i)}$ which extends the union $\bigcup\limits_{j
< i}f_j$.
\end{enumerate}
there exists a function $f : \lambda \rightarrow \kappa$ such that
for all $\delta \in S$, $f\restriction_{\rg(\eta_{\delta})} \in
P_{\delta}$.
\end{definition}

\begin{proposition}
\label{stronguni} Let $\lambda=\gamma^+$, $\gamma$ regular and let
$\bar{\eta}=\left< \eta_{\delta} : \delta \in S \right>$ be a tree-like
ladder system on $S$ such that $\bar{\eta}$ has $\gamma$-uniformization and
$\Diamond_{\gamma}$ holds. Then $\bar{\eta}$ has strong
$\gamma$-uniformization.
\end{proposition}

\begin{proof} Let $\bar{P}$ be given as stated and let $J$ be stationary
in $\gamma$ such that $\Diamond_{\gamma}(J)$ holds. For
simplicity we shall identify $J$ with $\gamma$. Thus there exists
a system of diamond functions $\bar{h}=\left< h_{\delta} : \delta
\rightarrow \gamma : \delta < \gamma \right>$ such that for
every function $h: \gamma \rightarrow \gamma$ the set $\{ \delta
< \gamma : {h\restriction_{\delta}}=h_{\delta} \}$ is
stationary in $\gamma$. For each $\delta \in S$ we define for $i
< \gamma$
\[ h_i^{\delta} : \rg(\eta_{\delta}\restriction_i) \rightarrow
\gamma \text{ via } \eta_{\delta}(j) \mapsto h_i(j). \] Moreover,
we put for $\delta \in S$, $E_{\delta}=\{ i < \gamma :
h_i^{\delta} \subseteq f \text{ for some } f \in P_{\delta} \}$
and let $g_i^{\delta} \in P_{\delta}$ be such that
$g_i^{\delta}\restriction_{\rg(\eta_{\delta}\restriction_i)}=h_i^{\delta}$
for $i \in E_{\delta}$. Note that $E_{\delta} \not= \emptyset$
since $\Diamond_{\gamma}$ holds. Define $f_{\delta} :
\rg(\eta_{\delta}) \rightarrow H(\gamma)$ for $\delta \in S$ as
follows
\[ f_{\delta}(\eta_{\delta}(i))=\left<
g_j^{\delta}\restriction_{\rg(\eta_{\delta}\restriction_{(i+1)})}
: j \leq i, j \in E_{\delta} \right>. \] Here $H(\gamma)$
denotes the class of sets hereditarily of cardinality $< \gamma$.
Note that $H(\gamma)$ has size $\leq \gamma$. By the
$\gamma$-uniformization of $\bar{\eta}$ we can find $F: \lambda
\rightarrow H(\gamma)$ such that for all $\delta \in S$ there
exists $\alpha_{\delta}< \gamma$ such that for all
$\alpha_{\delta} \leq i < \gamma$ we have
$f_{\delta}(\eta_{\delta}(i))=F(\eta_{\delta}(i))$. For $i <
\gamma$ let
\[ F(\eta_{\delta}(i))=\left< G_j^{\eta_{\delta}(i)} : j \leq i, j
\in E_{\eta_{\delta}(i)} \right> \] for some
$E_{\eta_{\delta}(i)} \subseteq \gamma$. Note, that
$F(\eta_{\delta}(i))$ does not depend on $\delta$ and that $F$ is
well-defined since the ladder system is tree-like.\\ We now define
$f: \lambda \rightarrow \gamma$ by defining $f:
\rg(\eta_{\delta}) \rightarrow \gamma$ for every $\delta \in S$.
Clearly this is enough. In order to define $f$ on
$\rg(\eta_{\delta})$ for fixed $\delta \in S$ we use induction on
$i < \gamma$. For $i=0$ choose any member $u \in
P_{\eta_{\delta}(0)}$ and put
$f(\eta_{\delta}(0))=u(\eta_{\delta}(0))$. Now assume that
$f(\eta_{\delta}(j))$ has been defined for $j < i$ and $\delta
\in S$ such that for $j < i$,
$f\restriction_{\rg(\eta_{\delta}\restriction_{(j+1)})} \in
P_{\eta_{\delta}(j)}$. Put $\bar{f}_{\delta}=\{
f(\eta_{\delta}(j)) : j < i \}$ and let \[ J_i^{\delta}=\{ j
\in E_{\eta_{\delta}(i)} : G_j^{\eta_{\delta}(i)}
\restriction_{(\rg(\eta_{\delta}\restriction_i))} \subseteq
\bar{f}_{\delta} \}. \] If $j_i^{\delta}=\min(J_i^{\delta})$
exists then let
$f(\eta_{\delta}(i))=G_{j_i^{\delta}}^{\eta_{\delta}(i)}(\eta_{\delta}(i))$.
If $\min(J_i^{\delta})$ does not exist, then we distinguish
between two cases: if $i$ is a limit ordinal, then (iii) from
Definition \ref{defuni} implies that there is $f_i\in
P_{\eta_{\delta}(i)}$ which extends $\bigcup\limits_{j <
i}f\restriction_{\rg(\eta_{\delta}\restriction_{(j+1)})}$. If $i$
is a successor ordinal, then (ii) from Definition \ref{defuni}
ensures that there is $f_i \in P_{\eta_{\delta}(i)}$ extending
$f\restriction_{\rg(\eta_{\delta}\restriction_i)}$. In both cases put
$f(\eta_{\delta}(i))=f_i(\eta_{\delta}(i))$. Note that $f$ is
well-defined since the ladder system is tree-like, hence
$\min(J_i^{\delta})$ exists if and only if $\min(J_j^{\nu})$ exists
for $\eta_{\delta}(i)=\eta_{\nu}(j)$ $(\delta, \nu \in S$, $i,j <
\gamma)$. It remains to check that
$f\restriction_{\rg(\eta_{\delta})} \in P_{\delta}$. By the
uniformization we have for $\delta \in S$ that for $i \geq
\alpha_{\delta}$,
$G_j^{\eta_{\delta}(i)}=g_j^{\delta}\restriction_{\rg(\eta_{\delta}\restriction_{(i+1)})}$
for all $j \leq i$, $j \in E_{\eta_{\delta}(i)}$. We define
\[ h: \gamma \rightarrow \gamma \text{ via } h(j)=
f(\eta_{\delta}(j)). \] By $\Diamond_{\gamma}$ there exists
$\beta_{\delta} \geq \alpha_{\delta}$ such that
$h\restriction_{\beta_{\delta}}=h_{\beta_{\delta}}$ and hence
\[
f\restriction_{\rg(\eta_{\delta}\restriction_{\beta_{\delta}})}=h_{\beta_{\delta}}^{\delta}\restriction_{\rg(\eta_{\delta}\restriction_{\beta_{\delta}})}=h_{\beta_{\delta}}^{\delta}.
\]
Thus in the construction $j_{\beta_{\delta}}^{\delta}$ existed
since for example $\beta_{\delta} \in
J_{\beta_{\delta}}^{\delta}$. Therefore, by definition of $f$, we have
\[
f\restriction_{\rg(\eta_{\delta}\restriction_{(\beta_{\delta}+1)})}=g_{\beta_{\delta}}^{\delta}\restriction_{\rg(\eta_{\delta}\restriction_{(\beta_{\delta}+1)})}.
\]
By induction on $i \geq \beta_{\delta}$ we can now show that
$f\restriction_{\rg(\eta_{\delta})}=g_{\beta_{\delta}}^{\delta}
\in P_{\delta}$ and this finishes the proof. \end{proof}

We can now improve Theorem \ref{uniformization}.

\begin{theorem}
\label{stronguniext}
Let $T$ be a group of cardinality less
than $\lambda$ and $G$ a group which is
$(T,\lambda,\gamma)$-suitable. Moreover, assume that $\lambda=\gamma^+$, $S$ is
non-reflecting and $\bar{\eta}$ is a tree-like ladder system on $S$ that
has $\gamma$-uniformization. If $\Diamond_{\gamma}$ holds,  then there
exists a torsion-free group $H$ of size $\lambda$ such that
\begin{enumerate}
\item $H$ has a $\lambda$-filtration $\left< \bar{H}_\alpha : \alpha < \lambda
\right>$;
\item if $\alpha \in S$, then $\bar{H}_{\alpha + 1}/\bar{H}_{\alpha} \cong G$;
\item if $\alpha \not\in S$, then $\bar{H}_{\beta}/\bar{H}_{\alpha}$ is free
for all $\alpha \leq \beta$;
\item $\Ext(H,T)=0$.
\end{enumerate}
In fact, $H$ satisfies $\Ext(H,W)=0$ for all groups $W$ of
cardinality less than or equal to  $| T |$.
\end{theorem}

\begin{proof}
The proof is almost identical with the proof of \cite[Proposition 1.8]{ES} but for the convinience of the reader we state it briefly pointing out the major changes. Let $S$, $\bar{\eta}$ and $G$ be as stated and choose the
group $H_{\bar{\eta}}=\bigoplus\{ H_{\eta} : \eta \in B_{\bar{\eta}}\}/K$ as constructed above. As in the proof of Theorem \ref{uniformization} $H_{\bar{\eta}}$ is an almost-free non-free torsion-free group of size $\lambda$ which has the desired $\lambda$-filtration. It remains to show that $\Ext(H_{\bar{\eta}}, K)=0$ for all groups $K$ of size less than or equal to the size of $T$. Let $K$ be such a group and choose a short exact sequence
\begin{equation}
\tag{$*$}
0 \rightarrow K \overset{id}{\rightarrow} N \overset{\pi}{\rightarrow} H_{\bar{\eta}} \rightarrow 0.
\end{equation}
We have to show that $(*)$ splits, i.e. we have to find a splitting map $\varphi: H_{\bar{\eta}} \rightarrow N$ such that $\pi \circ \varphi=id\restriction_{H_{\bar{\eta}}}$. Choose any set function $u: H_{\bar{\eta}} \rightarrow N$ such that $\pi \circ u=id\restriction_{H_{\bar{\eta}}}$. As in the proof of \cite[Proposition 1.8]{ES} the splitting maps $\varphi$ of $\pi$ are in one-one correspondence with set mappings $h: H_{\bar{\eta}} \rightarrow N$ such that $h(0)=0$ and for all $x,y \in H_{\bar{\eta}}$ and $z \in \Z$
\begin{enumerate}
\item $zh(x)-h(zx)=zu(x)-u(zx)$ and
\item $h(x) + h(y) - h(x+y)=u(x)+u(y)-u(x+y)$
\end{enumerate}
holds. For a subgroup $H$ of $H_{\bar{\eta}}$ we denote by $\Trans(H,K)$ the set of all set mappings $h: H \rightarrow N$ satisfying conditions (i) and (ii) from above for all $x,y \in H$ and $z \in \Z$. Thus $(*)$ splits if and only if $\Trans(H_{\bar{\eta}},K)$ is non-empty.\\
It is now easy to see that we can identify the elements of $\Trans(H_{\eta},K)$ with functions from $F_{\l(\eta)}$ to $K$ for each $\eta \in B_{\bar{\eta}}$. Remember that $\left< F_{\alpha} : \alpha \leq \gamma \right>$ was the chain we fixed as in Definition \ref{defsuitable} for $G$. For $\delta \in S$, $i < \gamma$ and $h \in \Trans(H_{\eta_{\delta}\restriction_i}, K)$ let $\seq(h)$ be defined as follows:
\[ \seq(h) : \rg(\eta_{\delta}\restriction_i) \rightarrow K \text{ via } \eta_{\delta}(j) \mapsto h\restriction_{F_{\l(\eta_{\delta}\restriction_j})} \text{ } (j < i). \]
For $\delta \in S$ let $P_{\delta}=\{ \seq(h) : h \in \Trans(H_{\eta_{\delta}},K) \}$ and for $i < \gamma$ put $P_{\eta_{\delta}(i)}=\{ \seq(h) : h \in \Trans(H_{\eta_{\delta}\restriction_{(i+1)}},K) \}$. Let $P_{\alpha}$ be trivial if it has not been defined yet ($\alpha < \lambda$). By Proposition \ref{stronguni} the ladder system $\bar{\eta}$ has strong $\gamma$-uniformization and it is easy to check that the system $\bar{P}=\left< P_{\alpha} : \alpha \in \lambda \right>$ satisfies the conditions of Definition \ref{defuni} since $F_n$ and $F_n/F_m$ are free for $m<n\leq\gamma$. Thus there exists a function $f: \lambda \rightarrow K$ such that for all $\delta \in S$, $f\restriction_{\rg(\eta_{\delta})} \in P_{\delta}$ since $K$ is of size less than or equal to $\gamma$. We now define $h: H_{\bar{\eta}} \rightarrow K$ by putting $h\restriction_{H_{\eta_{\delta}}}=f\restriction_{\rg(\eta_{\delta})}$ and clearly $h$ is well-defined and belongs to $\Trans(H_{\bar{\eta}},K)$ and therefore $(*)$ splits.
\end{proof}

\section{The Forcing Theorem}

Before we state the main theorem of this section let us describe
our strategy in order to make the statement of the main theorem
plausible. Using class forcing we will construct a model of $ZFC$
satisfying $GCH$ in which for every regular cardinal $\lambda$
there exists a sequence of stationary non-reflecting subsets
$S_{\alpha}$ of $\lambda$ of length $\lambda^+$ on which we have
"enough" uniformization for some ladder system. Using this and the
existence of $(T,\lambda,\gamma)$-suitable groups (for some
particular $\gamma$) we can then construct, for a given torsion
group $T$, a sequence of torsion-free groups $G_{\alpha}$ ($\alpha
< \lambda^+$) of cardinality $\lambda$ satisfying
$\Ext(G_{\alpha},T)=0$. These $G_{\alpha}$ will have
$\lambda$-filtrations $\left< G_{\alpha,\delta}: \delta < \lambda
\right>$ whose successive quotients satisfy $\Ext(G_{\alpha,
\delta+1}/G_{\alpha,\delta},T)\not=0$ for $\delta \in S_{\alpha}$. However, we will be able to
show that all these groups together do not fit into a single
group $G \in \tf_{\lambda}(T)$ via embedding since this would force $\Ext(G,T)\not=0$.
Thus there can not be any $\lambda$-universal group for $T$.

\begin{theorem}
\label{forcing} Let $V$ be a model of $ZFC$ in which the
generalized continuum hypothesis $GCH$ holds. Then for some class
forcing $\P$ not collapsing cardinals and preserving $GCH$ the
following is true
in $V^{\P}$:\\

If $\lambda > \gamma$ are infinite regular cardinals such that
$\gamma=\cf(\mu)$ if $\lambda=\mu^+$ then
\begin{enumerate}
\item there is a normal ideal $J=J_{\gamma}^{\lambda}$ on
$\lambda$;
\item there is a stationary subset $S=S_{\gamma}^{\lambda}$ of
$\lambda$ such that $S \in J$;
\item if $\delta \in S$, then $\cf(\delta)=\gamma$;
\item $S$ is non-reflecting, i.e. $S \cap \alpha$ is not
stationary in $\alpha$ for every $\alpha < \lambda$;
\item if $S^{\prime} \subseteq S$ is stationary in $\lambda$, then
there is a stationary $S^* \in J$ such that $S^*\subseteq
S^{\prime}$;
\item if $S^{\prime} \subseteq S$ and $S^{\prime} \not\in J$, then
$\Diamond_{S^{\prime}}$ holds;
\item if $S^{\prime} \subseteq S$ is stationary and $S^{\prime} \in
J$, then there exists a tree-like ladder system on $S^{\prime}$
which has $\mu$-uniformization for all $\mu < \kappa$ if
$\lambda=\kappa^+$ and $\kappa$ is singular, and for all $\mu <
\lambda$ otherwise;
\item there are $S_{\epsilon}=S_{\gamma,\epsilon}^{\lambda} \in J$ for
$\epsilon < \lambda^+$ such that
\begin{enumerate}
\item if $\eta < \epsilon< \lambda^+$, then $S_{\epsilon}
\backslash S_{\eta}$ is bounded;
\item if $\epsilon < \lambda^+$, then $S_{\epsilon+1} \backslash
S_{\epsilon}$ is stationary;
\item $J=\{ S^{\prime} \subseteq S : \exists \epsilon < \lambda^+
\forall \epsilon < \nu < \lambda^+$ $S^{\prime} \backslash
S_{\nu}$ is not stationary $\}$.
\end{enumerate}
\end{enumerate}
Moreover, if $\lambda=\cf(\lambda)> \gamma$, then there is a
stationary $S^*=S^{\lambda,*}_{\gamma}$ such that
\renewcommand{\labelenumi}{{\rm (}\arabic{enumi}{\hspace*{0.1cm}\rm)}}
\begin{enumerate}
\item if $\alpha \in S^*$, then $\cf(\alpha)=\gamma$;
\item $S^*$ is non-reflecting;
\item $\Diamond_{S^*}$ holds.
\end{enumerate}
\end{theorem}

The proof of Theorem \ref{forcing} will be divided into several
steps. First we deal with each regular $\lambda$ seperately and
then use Easton-support iteration to put the forcings together.
We will assume a knwoledge of forcing and our notation follows
that of \cite{J} with the exception that $p \leq q$ means that
the condition $q $ is stronger than the condition $p$. Let
$\lambda$ be a cardinal. Recall that a poset $P$ is called {\it
$\lambda$-complete} if for every $\kappa < \lambda$, every
ascending chain \[ p_0 \leq p_1 \leq \cdots \leq p_{\alpha} \quad
(\alpha < \kappa)\] has an upper bound. Moreover, $P$ is said to
be {\it $\lambda$-strategically complete} if Player I has a
winning strategy in the following game of length $\kappa$ for
every $\kappa < \lambda$. Players I and II alternately choose an
ascending sequence \[ p_0 \leq p_1 \leq \cdots \leq p_{\alpha}
\quad (\alpha < \kappa) \] of elements of $P$, where Player I
chooses at the even ordinals; Player I wins if and only if at each
stage there is a legal move and the whole sequence, $\left<
p_{\alpha} : \alpha < \kappa \right>$ has an upper bound (see
also \cite[Definition A1.1]{S1}). Note that, if $P$ is
$\lambda$-strategically complete and $G$ is generic over $P$, then
$V[G]$ has no new functions from $\kappa$ into $V$ for all
$\kappa < \lambda$, hence cardinals $\leq \lambda$ and their
cofinalities are preserved.

\begin{proposition}
\label{costationary} Let $\lambda$ be a regular cardinal and
assume $\lambda^{< \lambda}=\lambda$. For any regular $\kappa <
\lambda$, there exists a poset $\Q$ of cardinality $\leq \lambda$
which is $\lambda$-strategically complete (and hence preserves
all cardinals and preserves cofinalities $\leq \lambda$) and is
such that, for $G$ generic over $P$, in $V[G]$ there exists a
non-reflecting stationary and co-stationary subset $S$ of
$\lambda$ such that every member of $E$ has cofinality $\kappa$.
(Here, co-stationary means that the set $\lambda \backslash S$ is
also stationary).
\end{proposition}

\begin{proof}
The proof is similar to the proof of \cite[Lemma 2.3]{ES} but for
the convenience of the reader we state it briefly. We let $Q$ be
the set of all functions $q:\alpha \rightarrow 2=\{ 0,1\}$
($\alpha < \lambda$) such that $q(\mu)=1$ implies that
$\cf(\mu)=\kappa$ and such that for all limits $\delta \leq
\alpha$, the intersection of $q^{-1}[1]$ with $\delta$ is not
stationary in $\delta$. Then, for $G$ generic over $P$,
\[ S=\bigcup\{ q^{-1}[1] : q \in G \} \]
will be the desired set. We have to prove that $S$ is stationary
and co-stationary in $\lambda$. Hence, assume that $q$ forces $f$
is the name of a continuous increasing function $\bar{f}:\lambda
\rightarrow \lambda$; choose an ascending chain
\[ q_o \leq q_1 \leq \cdots \leq q_{\alpha} \quad (\alpha <
\kappa)\] such that for each $\alpha$ there exist
$\beta_{\alpha}, \gamma_{\alpha}$ such that $q_{\alpha} \Vdash
\bar{f}(\beta_{\alpha})=\gamma_{\alpha}$ and
\[ \dom(q_{\alpha}) \geq \gamma_{\alpha} > \dom(q_{\mu}) \]
for all $\mu < \alpha$. Let $\delta=\sup\{ \gamma_{\alpha} :
\alpha < \kappa\}=\sup \{\dom(q_{\alpha}) : \alpha < \kappa\}$
and let
\[ q_i=\bigcup\{ q_{\alpha} : \alpha < \kappa\} \cup
\{(\delta,i)\},\] for $i=0,1$. Then $q_i \in Q$ $(i=0,1$) since
$q_i^{-1}[1]$ is not stationary in $\delta$, because $\delta$ has
cofinality
$\kappa$. Moreover, $q_1 \Vdash \delta \in \rg(f) \cap S$ and $q_0 \Vdash \delta \in \rg(f) \cap (\lambda \backslash S)$.\\
Since $Q$ has cardinality $\leq \lambda$, it preserves cardinals
$> \lambda$. To show that all cardinals $\leq \lambda$ are
preserved (and their cofinalities), it suffices to prove that $Q$
is $\lambda$-strategically complete. Let $\tau < \lambda$ be a
limit ordinal. Let Player I choose $q_{\alpha}$ for even $\alpha$
such that $\dom(q_{\alpha})$ is a successor ordinal, say
$\delta_{\alpha} + 1$, and $q_{\alpha}(\delta_{\alpha})=0$.
Moreover, at limit ordinals $\alpha$ he chooses $q_{\alpha}$ to
have domain $=\sup\{ \delta_{\beta} : \beta < \alpha \}+1$. Then
$q=\bigcup\{q_{\alpha} : \alpha < \mu\}$ is a member of $Q$
because $\{ \delta_{\alpha} : \alpha < \mu, \alpha \text{ even
}\}$ is a cub in $\dom(q)$ which misses $q^{-1}[1]$. This is a
winning strategy for Player I and thus $Q$ is
$\lambda$-strategically complete.
\end{proof}

The next proposition is a collection of results from \cite{S4},
\cite{S1} and \cite{S2} (see also \cite{S5}).

\begin{proposition}
\label{oldforcing} Let $\lambda > \gamma$ be regular infinite
cardinals. Moreover, assume $\lambda^{< \lambda}=\lambda$, $2^{\lambda}=\lambda^+$ and let
$S$ be a non-reflecting, stationary and co-stationary subset of $\lambda$ such that
each member of $S$ has cofinality $\gamma$. Furthermore, let
$\gamma = \cf(\kappa)$ if $\lambda=\kappa^+$. Then there exists a
poset $P$ of cardinality $\leq \lambda^+$ which is
$\lambda$-strategically complete, satisfies the $\lambda^+$ chain
condition, adds no new sequences of length $< \lambda$ and has the
following properties:
\begin{enumerate}
\item $S$ is non-reflecting, stationary and co-stationary in $\lambda$ in $V^P$;
\item if $\lambda$ is inaccessible, then every ladder system on
$S$ has $\mu$-uniformization for all $\mu < \lambda$; in particular, there exists a tree-like ladder system on $S$;
\item if $\aleph_2 \leq \lambda=\kappa^+$ and $\kappa$ is regular,
then every ladder system on $S$ has $\mu$-uniformization for all
$\mu < \lambda$; in particular, there exists a tree-like ladder system on $S$;
\item if $\lambda=\aleph_1$, then there is a tree-like ladder system on
$S$ which has $\mu$-uniformization for all $\mu < \lambda$;
\item if $\lambda=\kappa^+$ and $\kappa$ is singular, then there
is a tree-like ladder system on $S$ which has $\mu$-uniformization
for all $\mu < \kappa$.
\end{enumerate}
\end{proposition}

\begin{proof}
For $\lambda$ inaccessible the proof is contained in \cite[Case A]
{S1} and also for the case of $\lambda=\kappa^+$, $\kappa$
regular (see \cite[Case B]{S1}). For $\lambda=\aleph_1$ see
\cite[V 1.7]{S4} and for $\lambda=\kappa^+$, $\kappa$ singular
see \cite[2.10, 2.12]{S2}. Moreover, simpler versions with less complicated and comprehensive proofs can be
found in \cite{S3} for all cases if we drop the requirements "for
every ladder system..." which is in fact not really needed for our purposes. Finally, let us remark that the
co-stationarity is only needed for $\lambda$ being the successor
of a regular cardinal or inaccessible.
\end{proof}

\begin{theorem}
\label{singleforcing} Let $\lambda$ be a regular cardinal such
that $\lambda^{< \lambda}=\lambda$ and $2^{\lambda}=\lambda^+$. Then there is a poset $P$ of
cardinality $\leq \lambda^+$ satisfying the $\lambda^+$-chain
condition which is $\lambda$-strategically complete and adds no
new sequences of length $< \lambda$ such that in $V^P$ for every
regular $\gamma < \lambda$ with $\gamma =\cf(\kappa)$ if
$\lambda=\kappa^+$ the statements (i) to (viii) and (1) to (3) of
Theorem \ref{forcing} hold.
\end{theorem}

Before we prove Theorem \ref{singleforcing} let us show why this
is enough to prove Theorem \ref{forcing}.

\begin{proof}(of Theorem \ref{forcing})\\
We start with a model $V$ of $ZFC$ satisfying the generalized
continuum hypothesis $GCH$. For any ordinal $\alpha$ let
$P_{\alpha}=\left< P_j,\dQ_i : j \leq \alpha, i < \alpha \right>$
be an iteration with Easton support; i.e. we take direct limits
when $\aleph_{\alpha}$ is regular and inverse limits elsewhere or
equivalently we have bounded support below inaccessibles and full
support below non-inaccessibles. For any ordinal $i$, let $\dQ_i$
be the forcing notion in $V^{P_i}$ described in Theorem
\ref{singleforcing} for $\lambda=\aleph_i$ if $\aleph_i$ is
regular and let it be $0$ elsewhere. Let $P$ be the direct limit
of the $P_{\alpha}$ ($ \alpha \in \Ord$). We claim that $P$ has
the desired properties. The proof is very similar
to the proof of \cite[Theorem 2.1]{ES} and hence we will only
state the main ingredients which are needed.
\begin{enumerate}
\item For every $\kappa$ and Easton support iteration
$\left< P_j, \dQ_i : \kappa \leq j \leq \alpha, \kappa\leq i <
\alpha \right>$, if each $\dQ_i$ is $\kappa$-strategically
complete, then so is $P_{\alpha}$.
\item $P=P_{\alpha}*P_{\geq \alpha}$, where, in $V^{P_{\alpha}}$, $P_{\geq \alpha}$
is the direct limit of $P_{\beta}^{\alpha}$ ($\beta \in \Ord$),
with $P_{\beta}^{\alpha}$ the Easton support iteration $\left<
P_j^{\alpha},\dQ_i^{\alpha} : j \leq \beta, i < \beta \right>$
where $\dQ_i^{\alpha}=\dQ_{\alpha+i}$.
\item $|P_n|=1$ (for $n \in \omega$); if $\aleph_{\delta}$ is
singular, $|P_{\delta}| \leq \aleph_{\delta}$ and
$|P_{\delta}|\leq \aleph_{\delta +1}$ if $\aleph_{\delta}$ is
regular, hence inaccessible.
\item $P_{\geq \alpha}$ is $\aleph_{\alpha}$-strategically
complete, and $P_{\geq \alpha +1}$ is even $\aleph_{\alpha + n}$-
strategically complete for all $n \in \omega$.
\end{enumerate}
By construction (i) to (viii) and (1) to (3) of Theorem
\ref{forcing} are now satisfied in $V^P$. Note, that stationarity
is preserved in the iteration because $P_{\geq \alpha}$ is
$\aleph_{\alpha}$-strategically complete. It remains to prove
that $V^P$ is a model of $ZFC$ satisfying $GCH$ and preserving
cofinalities (and hence cardinals). This follows very similar as
in the proof of \cite[Theorem 2.1]{ES} and hence we will omit the
proof here and leave it as an exercise to the reader.\end{proof}

It remains to prove Theorem \ref{singleforcing}.

\begin{proof} (of Theorem \ref{singleforcing})
The proof follows from the results in \cite{S1} and \cite{S2} but for the convenience of the reader we shall give some details. If $\lambda=\kappa^+$ is a successor cardinal, then we are easily done since there is only one $\gamma$, namely $\gamma=\cf(\kappa)$ under consideration. We choose $P$ to be the two step iterated forcing of the two forcings from Proposition \ref{costationary} and from Proposition \ref{oldforcing} with $\gamma=\cf(\kappa)$. Moreover, we may assume that $P$ also forces the sets $S^*=S^{\lambda,*}_{\gamma}$ satisfying Theorem \ref{forcing} (1) to (3) by an initial forcing. Note that the assumptions on $\lambda$ in Theorem \ref{singleforcing} are satisfied by \cite[Ex 12, page 70]{HSW}. If $\lambda$ is inaccessible, then it is more complicated since we have to deal with all regular $\gamma < \lambda$. But this was already done in \cite[Case B]{S1} where a stronger version of Proposition \ref{oldforcing} was shown. It was proved that there is even a forcing notion $P$ such that for all regular $\gamma < \lambda$ and given non-reflecting stationary, co-stationary subsets $S_{\gamma}$ of $\lambda$ consisting of ordinals $\alpha \in S_{\gamma}$ of cofinality $\gamma$, every ladder system on $S_{\gamma}$ has $\mu$-uniformization for all $\mu < \lambda$. Using this stronger result and again forcing the sets $S_{\gamma}^{\lambda,*}$ satisfying Theorem \ref{forcing} (1) to (3) it remains to show that we can define the ideal $J=J^{\lambda}_{\gamma}$ satisfying Theorem \ref{forcing} (vi) and (viii) (point (vii) of Theorem \ref{forcing} is clear). \\
Our forcing $P$ (from \cite{S1} and \cite{S2}) is the result of a $(<\lambda)$-support iteration of length $\lambda^+$, say $\left< P_i, \dQ_j : i \leq \lambda^+, j < \lambda^+ \right>$. Let us assume that $P_{\gamma}$ forces the set $S=S_{\gamma}^{\lambda}$ and the tree-like ladder system $\bar{\eta}$ on it. In $V^{P_{\gamma}}$ there exists a sequence $\left< S_{\epsilon}=S^{\lambda}_{\gamma, \epsilon} : \epsilon < \lambda^+ \right>$ such that
\begin{enumerate}
\item $S_{\epsilon} \subseteq S$;
\item $\eta < \epsilon < \lambda^+$ implies $S_{\epsilon}\backslash S_{\eta}$ is bounded;
\item $\epsilon < \lambda^+$ implies $S_{\epsilon +1} \backslash S_{\epsilon}$ is staionary.
\end{enumerate}
Now we define $J=J_{\gamma}^{\lambda}$ as $J=\{ S^{\prime} \subseteq S : \exists \epsilon < \lambda^+
\forall \epsilon < \nu < \lambda^+$ $S^{\prime} \backslash
S_{\nu}$ is not stationary $\}$. For each $i < \lambda^+$, $\dQ_i$ forces $\mu_i$-uniformization for the ladder system $(\dA_i, \bar{\df}_i)$ where $\bar{\df}_i=\left< \df_{\delta}^i : \delta \in \dS_i \right>$. Here $\dS_i$ and $\df_{\delta}^i$ are $P_i$-names for a member of $^{\dA_{\delta}}(\mu_i)$. Thus we obtain $\Vdash$ `` $\dS_i \subseteq \dS$ is stationary and there is $\epsilon < \lambda^+$ such that $(\forall \epsilon < \eta < \lambda^+)(\dS_i \cap \dS_{\eta} \backslash \dS_{\epsilon}$ is not stationary ``. A condition in $\dQ_i$ is for instance given by $g: \alpha \rightarrow \mu_i$ such that $\delta \in \dS_i$, $\delta \leq \alpha$ implies $f_{\delta}^i \subseteq^* g$.\\ 
It remains to show that $\Diamond_{S^{\prime}}$ holds for $S^{\prime} \not\in J$. Choose $i < \lambda^+$ such that $S^{\prime}$ comes from $V^{P_i}$. For some $j \in (i,\lambda^+)$, $Q_j$ is adding $\lambda$ cohen reals and we can interprete it as adding a diamond sequence $\left< \rho_{\epsilon} : \epsilon \in S^{\prime} \right>$ by initial segments. Trivially, in $V^{P_{j+1}}$, $\Diamond_{S^{\prime}}$ holds and we may work in $V^{P_{j+1}}$ now. For $\chi$ large enough we can find for every $x \in H(\chi(\lambda))$ an increasing continuous sequence $\bar{N}=\left< N_i : i < \lambda \right>$ of elementary submodels of $H(\chi(\lambda), \epsilon, <^*)$ of cardinality less than $\lambda$ such that $x \in N_0$, $S^{\prime} \in N_0$ and $\bar{N}\restriction_{(i+1)} \in N_{i+1}$ for all $i < \lambda$. Let $E=\{ \delta < \lambda : N_{\delta} \cap \lambda = \delta \}$ which is a cub in $\lambda$. Thus, for $\delta \in E$, for every $p \in \left(P/P_{j+1}\right)\cap N_{\delta}$ there is a condition $p \leq q \in P/P_{j+1}$ which is $(N_{\delta}, P/P_{j+1})$ generic and forces a value to $G \cap N_{\delta}$. It is known that we can now replace the diamond sequence on $S^{\prime}$ which we have in $V^{P/P_{j+1}}$ by one that is preserved by forcing with $P/P_{j+1}$ since $P/P_{j+1}$ adds no new subsets of $\lambda$ of length less than $\lambda$ and by the strategically completeness. This finishes the proof.  
\end{proof}

\section{Application to Kulikov's question}
In this final section we show that the answer to Kulikov's
question is consistently no for large classes of groups (not
necessarily torsion) and cardinals $\lambda$.

\begin{definition}
Let $T$ be a group of cardinality less than $\lambda$ and
let $H$ be a torsion-free group of cardinality $\lambda$.
Moreover, let $\{H_\alpha : \alpha < \lambda \}$ be a
$\lambda$-filtration of $H$. Then
\[ S_{\gamma}^{\lambda}[H,T]=\{ \delta \in S : \Ext(H_{\delta+1}/H_{\delta},T)\not=0
\}. \]
\end{definition}

Let us remark that the definition obviously depends on the given
filtration. We could make it independent by defining an
equivalence relation and letting $\Gamma^S_H(T)=\{ E \subseteq S
: \text{ there exists a stationary set } R \subseteq S \text{ such
that } E \cap R=S_{\gamma}^{\lambda}[H,T] \cap R \}$. But since we
don't need this and since the filtration under consideration shall
always be clear from the context we avoid additional notations.

\begin{theorem}
\label{diamond} Let $T$ be a group of cardinality less
than $\lambda$ and let $H$ be a torsion-free group of size
$\lambda$. Let $S=S_{\gamma}^{\lambda}[H,T]$. If $\Diamond_{S}$
holds, then $\Ext(H,T)\not=0$.
\end{theorem}

\begin{proof} The proof is standard and can be found for example in
\cite[Theorem XII.1.15]{EM}. \end{proof}

We shall now work in the model $V^{\P}$ obtained in Theorem
\ref{forcing}. Thus all results shall be consistency results
with ZFC and GCH. The symbol $(V^{\P})$ indicates that the
statement holds in our model $V^{\P}$.

\begin{theorem}[$V^{\P}$]
\label{main} Assume that $T$ is a group of cardinality
less than $\lambda$ and $G$ is a torsion-free group which is
$(T,\lambda,\gamma)$-suitable with $\gamma < \lambda$ regular,
$\gamma =\cf(\mu)$ if $\lambda= \mu^+$. If $\lambda > 2^{(| T |^{\gamma})}$ or $\mu$ is regular, then there is no
$\lambda$-universal group for $T$.
\end{theorem}

\begin{proof} Assume that there is a $\lambda$-universal group $G^*$ for
$T$. If $\lambda > 2^{(| T |^{\gamma})}$, then for
$\epsilon < \lambda^+$ we apply Theorem \ref{uniformization} to
$S_{\gamma,\epsilon}^{\lambda}$, $T$, $G$ and the tree-like
ladder system $\bar{\eta}_{\epsilon}$ on
$S_{\gamma,\epsilon}^{\lambda}$ which exists by Theorem
\ref{forcing} (vii). Note that no $S_{\gamma,\epsilon}^{\lambda}$
reflects in any $\alpha < \lambda$ and that
$\bar{\eta}_{\epsilon}$ has $2^{(| T
|^{\gamma})}$-uniformization since $\lambda > 2^{(| T
|^{\gamma})}$ (and $\mu > 2^{(|T|^{\gamma})}$ if $\lambda=\mu^+$, $\mu$ singular). If $\mu$ is regular, then we may apply Theorem
\ref{stronguniext} instead of Theorem \ref{uniformization}. Note
that $\Diamond_{\mu}$ holds in $V^{\P}$ by Theorem \ref{forcing} (1) to (3). Hence for each $\epsilon
< \lambda^+$ we obtain a torsion-free group
$H_{\epsilon}=\bigcup\limits_{\alpha \in \lambda}H_{\epsilon,
\alpha}$ satisfying $\Ext(H_{\epsilon},T)=0$ and
$\Ext(H_{\epsilon, \alpha+1}/H_{\epsilon, \alpha},T)\not=0$ for
all $\alpha \in S_{\gamma}^{\lambda}[H_{\epsilon},T]=S_{\gamma,
\epsilon}^{\lambda}$ . By universality of $G^*$ there exist
embeddings $i_{\epsilon} : H_{\epsilon} \longrightarrow G^*$. We
claim that for each $\epsilon < \lambda^+$ we have $S_{\gamma,
\epsilon}^{\lambda}[H_{\epsilon},T] \subseteq
S_{\gamma}^{\lambda}[G^*,T]$ modulo a non-stationary set. To see
this choose a $\lambda$-filtration $G^*=\bigcup\limits_{\alpha
\in \lambda}G_{\alpha}$ of $G^*$ such that $\Ext(G_{\alpha
+1}/G_{\alpha}, T)=0$ if and only if for some $\beta > \alpha$ we
have $\Ext(G_{\beta}/G_{\alpha},T)=0$. Fix $\epsilon <
\lambda^+$, then there is a cub $C_{\epsilon} \subseteq \lambda$
such that for all $\alpha \in C_{\epsilon}$ we have $H_{\epsilon,
\alpha}=G_{\alpha}\cap H_{\epsilon}$. Thus for $\alpha < \beta
\in C_{\epsilon}$ it follows that $H_{\epsilon,
\beta}/H_{\epsilon, \alpha}=(G_{\beta}\cap
H_{\epsilon})/(G_{\alpha}\cap H_{\epsilon})\subseteq
G_{\beta}/G_{\alpha}$ and hence $\Ext(H_{\epsilon,
\beta}/H_{\epsilon, \alpha}, T)\not=0$ implies
$\Ext(G_{\beta}/G_{\alpha},T)\not=0$. Therefore also
$\Ext(G_{\alpha+1}/G_{\alpha}, T)\not=0$ and $C_{\epsilon}
\subseteq S_{\gamma}^{\lambda}[G^*,T]$. Thus by the definition of
the normal ideal $J$ (see Theorem \ref{forcing} point (viii)) we
have $\bar{S}=S_{\gamma}^{\lambda}[G^*,T] \not\in J$ and therefore
$\Diamond_{\bar{S}}$ holds by Theorem \ref{forcing} point (vi).
Hence $\Ext(G^*,T)\not=0$ by Theorem \ref{diamond} - a
contradiction. \end{proof}

Before we prove some corollaries let us remark that in $V^{\P}$ the general continuum hypothesis $GCH$ holds. Hence, for a group $T$ we have $2^{(| T |^{\gamma})}
\leq max\{{\gamma^+}^+,{| T |^+}^+\}$ which the reader should keep in mind.

\begin{corollary}[$V^{\P}$]
\label{limit} Let $T$ be a group of cardinality less
than $\lambda$ which is not cotorsion. If $\lambda$ is strongly
inaccessible, then there is no $\lambda$-universal group for $T$.
\end{corollary}

\begin{proof} Since $\lambda$ is strongly inaccessible it is a limit
ordinal and we may choose $\gamma=\omega$. Moreover, for every
$\alpha < \lambda$ we have $2^{\alpha} < \lambda$, hence $\lambda
> 2^{(| T |^{\omega})}$. Lemma \ref{notcotorsion} implies
that there is a $(T,\lambda,\omega)$-suitable group for $T$ and
hence Theorem \ref{main} shows that there is no
$\lambda$-universal group for $T$. \end{proof}

\begin{corollary}[$V^{\P}$]
\label{cofinomega} Let $T$ be a group such that $| T |^+<
\lambda$ which is not cotorsion. If $\lambda=\mu^+$ and
$\cf(\mu)=\omega$, then there is no $\lambda$-universal group for
$T$.
\end{corollary}

\begin{proof} Since $\cf(\mu)=\omega$, we have $\lambda > 2^{(| T
|^{\omega})}$ and hence we may choose $\gamma =\omega$ and
apply Theorem \ref{main} to see that there is no
$\lambda$-universal group for $T$. Note, that there is a
$(T,\lambda,\omega)$-suitable group by Lemma \ref{notcotorsion}.
\end{proof}

\begin{corollary}[$V^{\P}$]
Let $T$ be a group of cardinality less than $\lambda$ which is not cotorsion and let $\lambda >
2^{(| T |^{\omega})}$ be regular. Moreover, let $T$ have an
epimorphic image of size less than or equal to $\cf(\mu)$ if
$\lambda=\mu^+$ (e.g. $T$ torsion, $T$ mixed splitting or $| T
| \leq \cf(\mu)$). If $2^{(| T |^{\cf(\mu)})} < \lambda$
then there exists no $\lambda$-universal group for $T$.
\end{corollary}

\begin{proof} We choose $\gamma=\omega$ if $\lambda$ is a limit ordinal
and $\gamma=\cf(\mu)$ if $\lambda=\mu^+$. By Theorem \ref{forcing} (1) to (3)
there exists a stationary non-reflecting set $S \subseteq \gamma$
consisting of limit ordinals of cofinality $\omega$ such that
$\Diamond_S$ holds. Let $H$ be the epimorphic image of $T$ as
stated. By assumption we may apply Proposition \ref{diamondsuitable} to
$S$, $H$ and $\lambda, \gamma$ to obtain a
$(T,\lambda,\gamma)$-suitable group $G$. Since $\lambda
> 2^{(| T |^{\gamma})}$ we apply Theorem \ref{main} to see
that there is no $\lambda$-universal group for $T$. \end{proof}

\begin{corollary}
Let $T$ be a group which is not cotorsion and $\lambda$ a
cardinal. Then there exists a regular uncountable cardinal
$\delta > \lambda$ such that there exists no $\delta$-universal
group for $T$.
\end{corollary}

\begin{corollary}[$V^{\P}$]
Let $T$ be a torsion group of regular cardinality which is not cotorsion. If $\lambda > | T |$ is regular, then there is no $\lambda$-universal
group for $T$.
\end{corollary}

\begin{proof} Let $T$ and $\lambda$ be as stated. If $\lambda$ is strongly inaccessible , then Corollary \ref{limit} gives the claim. Hence assume that $\lambda=\kappa^+$. Thus $| T | \leq \kappa$ and by Proposition \ref{diamondsuitable} there exists a $(T,\lambda,\cf(\kappa))$-suitable group $G$ for $T$. If $\kappa$ is regular then Theorem \ref{main} shows that there is no $\lambda$-universal group for $T$. If $\kappa$ is singular, then $\kappa> (| T |^+)^+$ since $| T |$ is regular. Hence $\lambda > 2^{(| T |^{\cf(\kappa)})}$ and again Theorem \ref{main} gives the claim. \end{proof}

\goodbreak

\end{document}